\documentclass[12pt]{amsart}
\usepackage{amsmath}
\usepackage{amsfonts}
\usepackage{latexsym}
\usepackage{amssymb}
\usepackage{enumerate}
\usepackage[dvips]{graphics}

\newcommand{\tc}{{\rm\bf {TC}}}

\newcommand{\comment}[1]{}

\newtheorem{theorem}{Theorem}

\newtheorem{lemma}[theorem]{Lemma}

\newtheorem{corollary}[theorem]{Corollary}
\newtheorem{definition}{Definition}

\begin{document}
\title{Collision Free Motion Planning on Graphs}
\author[M. ~Farber]{M.~Farber}
\address{Department of Mathematics, University of Durham, Durham DH1 2LE, UK}
\email{Michael.Farber@durham.ac.uk}


\date{June 18, 2004}

\keywords{Topological robotics, Motion planning algorithm, Configuration spaces
of graphs}

\begin{abstract} A topological theory initiated in \cite{F1}, \cite{F2} uses methods of
algebraic topology to estimate numerically the character of instabilities
arising in motion planning algorithms. The present paper studies {\itshape
random} motion planning algorithms \index{random motion planning algorithm} and
reveals how the topology of the robot's configuration space influences their
structure. We prove that the topological complexity of motion planning $\tc(X)$
coincides with the minimal $n$ such that there exists an $n$-valued random
motion planning algorithm for the system; here $X$ denotes the configuration
space. We study in detail the problem of collision free motion of several
objects on a graph $\Gamma$. We describe an explicit motion planning algorithm
for this problem. We prove that if $\Gamma$ is a tree and if the number of
objects is large enough, then the topological complexity of this motion
planning problem equals $2m(\Gamma)+1$ where $m(\Gamma)$ is the number of the
essential vertices of $\Gamma$. It turns out (in contrast with the results on
the collision free control of many objects in {\itshape space} \cite{FY}) that
the topological complexity is independent of the number of particles.
\end{abstract}

\maketitle

%


\section{Introduction}

Algorithmic motion planning in robotics is a well established discipline which
provides (1) a wide variety of general-purpose efficient theoretical
algorithms, (2) more efficient algorithms designed for a number of special
situations and (3) some practical solutions which work reasonably well in
fairly involved scenarios. We refer to \cite{S} for a recent survey and to
\cite{L} for a comprehensive textbook.

In general, one is given a moving system $R$ with $k$ degrees of freedom and a
two or three-dimensional workspace $V$. The geometry of $R$ and of $V$ is given
in advance which determines the configuration space of the system, $X$. The
latter is a subset of $\mathbf R^k$ consisting of all placements (or
configurations) of the system $R$, each represented by a tuple of $k$ real
parameters, such that in this placement $R$ lies fully in $V$. For simplicity
we may restrict our attention to a single connected component of $R$, the one
containing a prescribed initial placement of $R$.

Being a subset of the Euclidean space $\mathbf R^k$, the configuration space
$X$ naturally inherits its {\itshape topology}. Many questions of control
theory depend solely on the configuration space $X$ viewed as a topological
space. One of the advantages of this approach is that different control
problems could be treated simultaneously for all systems having homeomorphic
configuration spaces. It is well-known that any real analytic manifold can be
realized as the configuration space of a simple mechanical system (linkage).
Therefore topological questions of robotics lead to interesting new topological
invariants of abstract manifolds.

We are interested in motion planning algorithms which work as follows: the
algorithm gets as its input the present and the desired states of the system
and it produces as the output a continuous motion of the system from its
current state to the desired state. It turns out that the topology of the
configuration space of the system imposes important restrictions on the
discontinuities of the robot motion as a function of the input data. We
emphasize that these are not discontinuities of the robot motion as a function
of time. The discontinuities which we study here are in the way the decision
(the whole motion) depends on the input data.

The approach to the motion planning problem described in \cite{F1}, \cite{F2}
was initially inspired by my discussions with Dan Halperin and Micha Sharir in
2000. It was also influenced by the well-known previous work of S. Smale
\cite{Sm} on algorithms of finding roots of algebraic equations.

Our results on topological robotics were published in mathematical journals
\cite{F1}, \cite{F2}, \cite{FTY} and in \cite{FY}. They led to an interesting
new topological invariant $\tc(X)$ of topological spaces. In robotics
applications, the number $\tc(X)$ has at least three different appearances.
Firstly, it is the minimal number of domains of continuity of any motion
planning algorithm for a system having $X$ as its configuration space.
Secondly, it is the minimal {\itshape order of instability}\index{order of
instability} (see \cite{F2}) which have motion planning algorithms in $X$. The
third interpretation (see \S \ref{sec:random} below) allows to measure $\tc(X)$
while relying on {\itshape random} motion planning algorithms: we show that
$\tc(X)$ is the minimal integer $n$ such that there exists an $n$ valued random
motion planning algorithm for the system.

The main part of the paper is devoted to a very specific motion planning
problem: simultaneous control of many objects whose motion is restricted by a
graph and the goal is to construct a motion planning algorithm avoiding
collisions between the objects. This problem was initially studied by R.
Ghrist, D. Koditschek and A. Abrams \cite{Gr}, \cite{GK}, \cite{Ab}. We
calculate here the topological complexity of the problem and describe an
explicit motion planning algorithm solving it.

For convenience of the reader we have included sections \S \S \ref{sec:one} -
\ref{sec:five} which give a brief description of some results of \cite{F1},
\cite{F2} used later in this paper.

\section{Motion Planning Algorithms}\label{sec:one}

Consider a mechanical system (robot) controlled by a motion planning algorithm. It is
supposed to function as follows: an operator introduces into the computer of the system
the current and the desired states of the system and the motion planning algorithm
determines a continuous motion of the system from its current state to the desired state.

Let $X$ be the configuration space of the system. We will always assume that
$X$ is path connected, i.e. any pair of points $A, B\in X$ may be joined by a
continuous path $\gamma$ in $X$. This means that it is possible to bring our
system, by a continuous movement, from any given configuration $A$ to any given
configuration $B$. This assumption does not represent a restriction since in
practical situations when the natural configuration space of a given system has
several connected components, we may simply restrict our attention to one of
them.

Given two points $A, B \in X$, one wants to connect them by a path in $X$; this
path represents a continuous motion of the system from one state to the other.
A motion planning algorithm is a rule (algorithm) that takes pairs of
configurations $(A,B) \in X \times X$ as an input and produces a continuous
path in $X$ from $A$ to $B$ as an output.

Let $PX$ denote the space of all continuous paths $\gamma: [0,1] \to X$,
equipped with the compact-open topology, and let $\pi: PX \to X \times X$ be
the map assigning the end points to a path: $\pi (\gamma) = (\gamma(0),
\gamma(1))$. Rephrasing the above discussion we see that a motion planning
algorithm is a {\itshape section} of this fibration
\begin{eqnarray}
s: X\times X\to PX, \quad \pi\circ s=1_{X\times X}.\label{section1}
\end{eqnarray}
Here $1_{X\times X}: X\times X\to X\times X$ denotes the identity map. The
algorithm associates the curve $s(A,B)(t)$, where $t\in [0,1]$, with any two
given configurations of the system $A, B\in X$.

Given a mechanical system, one asks {\itshape does there exist a continuous
motion planning algorithm for it}? In other words, whether it is possible to
find a continuous section (\ref{section1}). The answer is negative in most
cases as the following theorem proven in \cite{F1} states:

\begin{theorem}\label{contract} A globally defined continuous motion planning
algorithm (\ref{section1}) exists if and only if the configuration space $X$ of the
system is contractible.
\end{theorem}

\section{Topological Complexity of Motion Planning Algorithms}

In view of Theorem \ref{contract}, one expects that in general a motion
planning algorithm is only {\itshape piecewise continuous}. The following
definition describes an important class of piecewise continuous motion planning
algorithms having {\itshape finitely many} domains of continuity:

\begin{definition}\label{tame}
Let $X$ be a path-connected topological space. A motion planning algorithm
(\ref{section1}) is called {\itshape tame} if there exist finitely many subsets
\begin{eqnarray}
F_1, \dots, F_k\subset X\times X\label{decomposition}
\end{eqnarray}
such that the following conditions are
satisfied:
\begin{enumerate}
\item[(a)] the sets $F_1, \dots, F_k$ are pairwise disjoint $F_i\cap
F_j=\emptyset$, $i\not= j$, and cover $X\times X$, i.e. $ X\times X = F_1\cup
F_2\cup\dots \cup F_k;$ \item[(b)] each restriction $s|_{F_j}$ is continuous;
\item[(c)] each set $F_j$ is an ENR (see below).
\end{enumerate}
\end{definition}

Condition (a) means that the sets $F_1, \dots, F_k$ partition the total space
of all possible pairs $X\times X$. Condition (b) is the major continuity
assumption. Condition (c) is technical, it allows to avoid pathological
(exotic) decompositions. Recall, a topological space $Y$ is called an {\it
Euclidean Neighborhood Retract} (ENR) if it is homeomorphic to a subset of a
Euclidean space $Y'\subset {\mathbf R}^n$, such that $Y'$ is a retract of some
open neighborhood $Y'\subset U\subset {\mathbf R}^n$; in other words, $U\subset
{\mathbf R}^n$ is open and there exists a continuous map $r: U\to Y'$ such that
$r(y)=y$ for all $y\in Y'$. Such a continuous map $r$ is called a {\it
retraction}.

It is well-known that all manifolds and polyhedra are ENRs.

For a given algorithm (\ref{section1}) there may exist many different decompositions
(\ref{decomposition}) satisfying the conditions of Definition \ref{tame}.

\begin{definition}\label{algtopcompl}
{\itshape The topological complexity} of a tame motion planning algorithm $s: X\times
X\to PX$, $\pi\circ s=1_{X\times X}$ is defined as the minimal number $k$ of domains of
continuity $F_j$ which appear in Definition \ref{tame}.
\end{definition}

Given a concrete mechanical system one wishes to construct motion planning
algorithms for it with the minimal possible topological complexity. The problem
clearly depends only on the topology of the configuration space $X$ of the
system. This leads to the following purely topological notion:

\begin{definition}\label{maindef}
Let $X$ be a path-connected topological space. {\itshape The topological
complexity} of $X$ is defined as the minimal topological complexity of motion
planning algorithms in $X$.
\end{definition}

The topological complexity of a topological space $X$ coincides (for nice
spaces $X$) with the invariant $\tc(X)$ which was introduced in \cite{F1}. Its
definition (which is more convenient from the purely topological point of view)
appears in the following section.

\section{Topological Invariant $\tc(X)$}

\begin{definition}\label{opencover}
Let $X$ be a path-connected topological space. The number $\tc(X)$ is defined as the
minimal integer $r$ such that the Cartesian product $X\times X$ can be covered by $r$
open subsets
\begin{eqnarray}\label{cover} X\times X =
U_1\cup U_2\cup \dots \cup U_r,
\end{eqnarray} such that for any $i=1, 2,
\dots, r$ there exists a continuous map
\begin{eqnarray}
s_i: U_i\to PX\quad \mbox{with}\quad \pi\circ s_i =1_{U_i}.\label{local}
\end{eqnarray}
If no such $r$ exists, we set $\tc(X)=\infty$.
\end{definition}

It is shown in \cite{F1} that $\tc(X)$ depends only on the homotopy type of
$X$.

The invariant $\tc(X)$ admits an upper bound \cite{F1}:
\begin{eqnarray} \tc(X)\leq 2\dim X
+1.\label{upper}\end{eqnarray} The number $\tc(X)$ was computed in
\cite{F1},\cite{F2},\cite{FTY}, \cite{FY} for a number of important
configuration spaces appearing in robotics.

\begin{theorem}
Let $X$ be a connected smooth manifold. Then the topological complexity of $X$
(cf. Definition \ref{maindef}) coincided with $\tc(X)$ (cf. Definition
\ref{opencover}).
\end{theorem}

We refer to \cite{F2} for a proof.

A lower bound for $\tc(X)$ is based on the knowledge of the cohomology algebra
of $X$. To describe this result (which will be used later in this paper) we
first observe that the singular cohomology $H^\ast(X;{\mathbf R})=H^\ast(X)$ is
a graded ${\mathbf R}$-algebra with the multiplication
\begin{eqnarray}
\cup: H^\ast(X)\otimes H^\ast(X)\to H^\ast(X)\label{prod}
\end{eqnarray}
given by the cup-product, see \cite{DNF}, \cite{Sp}. The tensor product
$H^\ast(X)\otimes H^\ast(X)$ is again a graded ${\mathbf R}$-algebra with the
multiplication
\begin{eqnarray}\label{signs}
(u_1\otimes v_1)\cdot (u_2\otimes v_2) = (-1)^{|v_1|\cdot |u_2|}\,
u_1u_2\otimes v_1v_2.
\end{eqnarray}
Here $|v_1|$ and $|u_2|$ denote the degrees of cohomology classes $v_1$ and
$u_2$ correspondingly. The cup-product (\ref{prod}) is an algebra homomorphism.

\begin{definition} The kernel of homomorphism (\ref{prod}) is called
{\it the ideal of the zero-divisors} of $H^\ast(X)$. The {\it
zero-divisors-cup-length} of $H^\ast(X)$ is the length of the longest
nontrivial product in the ideal of the zero-divisors of $H^\ast(X)$.
\end{definition}

The next result is the main cohomological lower bound for the topological
complexity.

\begin{theorem}\label{lower} $\tc(X)$
is greater than the zero-divisors-cup-length of the cohomology algebra $H^\ast(X)$.
\end{theorem}

See \cite{F1} for a proof.

As an illustration we state the following result from \cite{F2} which is
relevant for the sequel:

\begin{theorem}\label{graph1}
Let $X$ be a connected graph. Then
\begin{eqnarray}\label{graph}
\tc(X) = \left\{
\begin{array}{lll}
1, & \mbox{if} & b_1(X)=0,\\
2, & \mbox{if} & b_1(X)=1,\\
3, & \mbox{if} & b_1(X)\geq 2.
\end{array}
\right.
\end{eqnarray}
Here $b_1(X)$ denotes the first Betti number of $X$.
\end{theorem}

\section{Order of Instability}\label{sec:five}

Besides the number of domains of continuity, the motion planning algorithms
could be characterized by their orders of instability.

Let $s:X\times X\to PX$, $\pi\circ s=1_{X\times X}$ be a tame motion planning
algorithm (cf. Definition \ref{tame}). Let $F_1, \dots, F_k\subset X\times X$
be pairwise disjoint subsets as in Definition \ref{tame}, i.e. such that
$s|_{F_j}$ is continuous, each $F_j$ is an ENR and union of the sets $F_j$
equals $X\times X$.

\begin{definition}
The {\itshape order of instability} of a motion planning algorithm $s$ is
defined as the smallest integer $r$ such that the subsets $F_1, \dots,
F_k\subset X\times X$ as above could be constructed in such a way that for any
sequence of $r+1$ indices $1\leq i_1< i_2<\dots <i_{r+1}\leq k$ one has
\begin{eqnarray}
\bar F_{i_1}\cap  \bar F_{i_2}\cap \dots \cap \bar F_{i_{r+1}}=\emptyset.
\end{eqnarray}
\end{definition}

Here $\bar F_j$ denotes the closure of $F_j$ in $X\times X$.

The order of instability represents a very important functional characteristic
of a motion planning algorithm. If the order of instability equals $r$ then
there exists a pair of initial - final configurations $(A,B)\in F_j$ such that
arbitrarily close to $(A,B)$ there exist $r-1$ pairs of configurations $(A_1,
B_1)$, $(A_2, B_2)$, $\dots$, $(A_{r-1}, B_{r-1})$ (which are all distinct
small perturbations of $(A,B)$), belonging to distinct sets $F_i$, $i\not=j$.
This means that {\itshape small perturbations of the input data $(A,B)$ may
lead to $r$ essentially distinct motions given by the motion planning algorithm
$s$.}

In practical situations one prefers to have motion planning algorithms with
order of instability as low as possible.

\begin{theorem}
Let $X$ be a connected smooth manifold. Then the order of instability of any
tame motion planning algorithm (\ref{section1}) is at least $\tc(X)$. Moreover,
there exists a motion planning algorithm with the order of instability
$\tc(X)$.
\end{theorem}

This theorem proven in \cite{F2} gives yet another way the topological
invariant $\tc(X)$ appears in robotics. In the next section we explain how the
topological quantity $\tc(X)$ shows in random algorithms.

\section{Random Motion Planning Algorithms}\label{sec:random}

According to Theorem \ref{contract}, {\itshape deterministic} continuous motion
planning algorithms $s:X\times X\to PX$ exist only when the configuration space
$X$ is contractible. Instead, one may work with {\itshape random} algorithms
solving the motion planning problem.

Let $X$ be a path-connected topological space. A {\itshape random $n$-valued
path}\index{random path} $\sigma$ in $X$ starting at $A\in X$ and ending at
$B\in X$ is given by an ordered sequence of paths $\gamma_1, \dots, \gamma_n\in
PX$ and an ordered sequence of real numbers $p_1, \dots, p_n\in [0,1]$ such
that each $\gamma_j:[0,1]\to X$ is a continuous path in $X$ starting at
$A=\gamma_j(0)$ and ending at $B=\gamma_j(1)$, and
\begin{eqnarray*}p_j\geq 0,
\quad p_1+p_2+\dots+p_n=1.
\end{eqnarray*}
One thinks of the paths $\gamma_1, \dots, \gamma_n$ as of the {\itshape states}
of $\sigma$ and of the number $p_j$ as being the {\itshape probability} that
the random path $\sigma$ is in state $\gamma_j$. Random path $\sigma$ as above
will be written as a formal linear combination
\begin{eqnarray}
\sigma= p_1\gamma_1+p_2\gamma_2 +\dots +p_n\gamma_n.\label{random}
\end{eqnarray}
Equality between $n$-valued random paths is understood as follows: the random
path (\ref{random}) is equivalent to $\sigma'= p'_1\gamma'_1+p'_2\gamma'_2
+\dots +p'_n\gamma'_n$ iff $p_j=p'_j$ for all $j=1, \dots, n$ and, besides,
$\gamma_j=\gamma'_j$ for all indices $j$ with $p_j\not=0$. In other words the
path $\gamma_j$ which appears with the zero probability $p_j=0$ could be
replaced by any other path starting at $A$ and ending at $B$.

We denote by $P_nX$ the space of all $n$-valued random paths in $X$. The space
$P_nX$ has a natural topology: it is a factor-space of a subspace of the
Cartesian product of $n$ copies of $PX\times [0,1]$. Note that the space $P_1X$
coincides with $PX$. The canonical map
\begin{eqnarray}
\pi: P_nX \to X\times X\label{canon}
\end{eqnarray}
assigns to a random path its initial and end points. Map (\ref{canon}) is
continuous.

\begin{definition}
An {\itshape $n$-valued random motion planning algorithm}\index{random motion
planning algorithm} is defined as a continuous map
\begin{eqnarray}
s: X\times X\to P_nX\label{randomalg} \end{eqnarray} such that $\pi\circ s=
1_{X\times X}$.
\end{definition}

Given a pair $(A,B)\in X\times X$, the output of the algorithm
\begin{eqnarray}
s(A,B) = p_1\gamma_1+ \dots+p_n\gamma_n\label{output}
\end{eqnarray}
is an ordered probability distribution on the paths between $A$ and $B$. In
other words, the algorithm $s$ produces the motion $\gamma_j$ with probability
$p_j$ where $j=1, \dots, n$.

\begin{theorem}\label{equiv}
Let $X$ be a path-connected metric space. Then the minimal integer $n$ such
that there exists an $n$-valued random motion planning algorithm $s:X\times
X\to P_nX$ in $X$ coincides with $\tc(X)$.
\end{theorem}
\begin{proof}
The following proof is an adjustment of the proof of Proposition 2 from
\cite{Sz}. Assume that there exists an $n$-valued random motion planning
algorithm $s:X\times X\to P_nX$ in $X$. The right hand side of formula
(\ref{output}) defines continuous real valued functions $p_j:X\times X\to
[0,1]$, where $j=1, \dots , n$. Let $U_j$ denote the open set
$p_j^{-1}(0,1]\subset X\times X$. The sets $U_1, \dots, U_n$ form an open
covering of $X\times X$. Setting $s_j(A,B) =\gamma_j$, one gets a continuous
map $s_j: U_j\to PX$ with $\pi\circ s_i=1_{U_j}$. Hence, $n\geq \tc(X)$
according to the definition of $\tc(X)$.

Conversely, setting $k=\tc(X)$, we obtain that there exists an open cover $U_1,
\dots, U_k\subset X\times X$ and a sequence of continuous maps $s_i:U_i\to PX$
where $\pi\circ s_i=1_{U_i}$, $i=1, \dots, k$. Extend $s_i$ to an arbitrary
(possibly discontinuous) mapping
\begin{eqnarray*}S_i:X\times X\to PX\end{eqnarray*} satisfying $\pi\circ
S_i=1_{X\times X}$. This can be done without any difficulty; it amounts in
making a choice of a connecting path for any pair of points $(A,B)\in X\times
X-U_i$. Next, one may find a continuous partition of unity subordinate to the
open cover $U_1, \dots, U_k$. It is a sequence of continuous functions $p_1,
\dots, p_k: X\times X\to [0,1]$ such that for any pair $(A,B)\in X\times X$ one
has
\begin{eqnarray*}
p_1(A,B) +p_2(A,B)+\dots + p_k(A,B) =1
\end{eqnarray*}
and the closure of the set $p_i^{-1}(0,1]$ is contained in $U_i$. We obtain a
continuous $k$-valued random motion planning algorithm $s: X\times X \to P_nX$
given by the following explicit formula
\begin{eqnarray}
s(A,B) =p_1(A,B)S_1(A,B)+\dots+ p_k(A,B)S_k(A,B).
\end{eqnarray}
The continuity of $s$ follows from the continuity of the maps $S_i$ restricted
to the domains $p_i^{-1}(0,1]$. This completes the proof.
\end{proof}

\section{Configuration Spaces of Graphs}

Let $\Gamma$ be a connected finite graph. The symbol $F(\Gamma, n)$ denotes the
configuration space of $n$ distinct particles on $\Gamma$. In other words,
$F(\Gamma, n)$ is the subset of the Cartesian product
\begin{eqnarray*}
\underbrace{\Gamma\times \Gamma\times \dots \times \Gamma}_{n \, \,
\mbox{\scriptsize times}} =\Gamma^{\times n}
\end{eqnarray*}
consisting of configurations $(x_1, x_2, \dots, x_n)$ where $x_i\in \Gamma$ and
$x_i\not= x_j$ for $i\not= j$. The topology of $F(\Gamma,n)$ is induced from
its embedding into $\Gamma^{\times n}$.

Configuration spaces of graphs were studied by R. Ghrist, D. Koditschek and A.
Abrams, see \cite{Gr}, \cite{GK}, \cite{Ab}.

To illustrate the importance of these configuration spaces for robotics one may
mention the control problems where a number of automated guided vehicles (AGV)
have to move along a network of floor wires \cite{GK}. The motion of the
vehicles must be safe: it should be organized so that the collisions do not
occur. If $n$ is the number of AGV then the natural configuration space of this
problem is $F(\Gamma, n)$ where $\Gamma$ is a graph. Here we idealize reality
by assuming that the vehicles have size 0 (i.e. they are points). Although this
assumption simplifies our discussion, it is in fact irrelevant for the
topological problems which we study.

The first question to ask is whether the configuration space $F(\Gamma, n)$ is connected.
Clearly $F(\Gamma, n)$ is disconnected if $\Gamma=[0,1]$ is a closed interval (and $n\geq
2$) or if $\Gamma=S^1$ is the circle and $n\geq 3$. These are the only examples of this
kind as the following simple lemma claims:

\begin{lemma}\label{connected}
Let $\Gamma$ be a connected finite graph having at least one essential vertex
 Then the configuration
space $F(\Gamma,n)$ is connected.
\end{lemma}

An {\itshape essential vertex}\index{essential vertex} is a vertex which is
incident to 3 or more edges.

One of the main results of topological robotics states that the configuration
spaces $F(\Gamma, n)$ are aspherical; see \cite{Gr}.

For the topological complexity of the configuration spaces one has:

\begin{theorem}\label{cor2}
Let $\Gamma$ be a connected graph having an essential vertex. Then the
topological complexity of $F(\Gamma, n)$ satisfies
\begin{eqnarray}\label{leq}
\tc(F(\Gamma, n)) \leq 2m(\Gamma)+1,
\end{eqnarray}
where $m(\Gamma)$ denotes the number of essential vertices in $\Gamma$.
\end{theorem}
\begin{proof}
By a theorem of R. Ghrist \cite{Gr}, the configuration space $F(\Gamma, n)$ has
homotopy type of a cell complex of dimension $\leq m(\Gamma)$. Combining this
result with the theorem about the homotopy invariance of $\tc(X)$ (see Theorem
3 in \cite{F1}) and the with the upper bound (\ref{upper}) one gets
(\ref{leq}).
\end{proof}

We prove below that equality holds in (\ref{leq}) in many cases.

We shall also see examples where (\ref{leq}) holds as a strict inequality.

\section{A Motion Planning Algorithm in $F(\Gamma,n)$}\label{algorithm}

In this section $\Gamma$ denotes a {\itshape tree}\footnote{Recall that tree is
a connected graph with no cycles.} having an essential vertex.
\begin{figure}[h]
\begin{center}
\resizebox{3.6cm}{4cm}{\includegraphics[167,380][503,761]{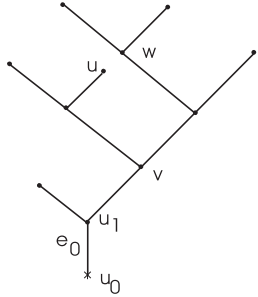}}
\end{center}
\caption{A partial order between the points of the tree.}\label{tree}
\end{figure}
Fix a univalent vertex $u_0\in \Gamma$ which will be called {\it the root}. Any
point in $\Gamma$ can be connected by a path to the root $u_0$ and this
connecting path is unique up to homotopy. The choice of the root determines a
partial order on $\Gamma$: we say that $x\succ y$, where $x, y\in \Gamma$ if
any path from $x$ to the root $u_0$ passes through $y$. Of course, $\succ$ is
only a partial order, i.e. there may exist pairs $x, y\in \Gamma$ such that
neither $x\succ y$, nor $y\succ x$. On Figure \ref{tree} we see $u\succ v$ and
$w\succ v$ however the points $u$ and $w$ are not comparable.

Let $e_0\subset \Gamma$ denote the root edge of $\Gamma$.

Fix a configuration $\alpha_0\in F(\Gamma, n)$ of $n$ distinct points lying on
$e_0$ and a continuous collision free motion connecting (in $F(\Gamma, n)$) any
pair of permutations of $\alpha_0$. Such motions exist as we assume that
$\Gamma$ has an essential vertex and hence the configuration space $F(\Gamma,
n)$ is connected (see Lemma \ref{connected}).

The algorithm works as follows. Let $\alpha=(A_1, A_2, \dots, A_n)\in F(\Gamma,
n)$ and $\beta=(B_1, B_2, \dots, B_n)\in F(\Gamma, n)$ be two given
configurations of $n$ distinct points on $\Gamma$. Let $A_{i_1}, \dots,
A_{i_r}$ be all the minimal elements (with respect to the order $\succ$) of the
set of points of $\alpha$. Here the notation is such that $i_1<i_2< \dots<i_r$.
First we move the point $A_{i_1}$ down to an interior point of the root edge
$e_0$. Next we move $A_{i_2}$ to the root edge $e_0$ and we continue moving
similarly the following points $A_{i_3}, \dots, A_{i_r}$.

As a result, after this first stage of the algorithm, all the minimal points of
$\alpha$ are transferred onto the root edge $e_0$. On the second stage we find
the minimal set among the remaining points of $\alpha$ and move them down, one
after another, to the edge $e_0$. Repeating in this way we find continuous
collision free motion of all the points of $\alpha$ moving them onto the
interior of the root edge $e_0$. We obtain a configuration of points
$\alpha'=(A_1', \dots, A'_n)$ which all lie in the interior of the root edge
$e_0$, in certain order.

For points lying on the root edge $e_0$ the partial order $\succ$ is a linear
order. Given a configuration $A=(A_1, \dots, A_n)$ where all points $A_i$ lie
on the edge $e_0$, there exists a unique permutation $\alpha: \{1, \dots,
n\}\to \{1, \dots, n\}$ such that $A_{\alpha_n}\succ A_{\alpha_{n-1}}\succ
\dots \succ A_{\alpha_1}.$ This permutation describes the order of the points
on the edge. It is obvious that any two configurations having the same
permutation can be connected by a continuous collision free path such that no
points leave $e_0$ in the process of motion.

Applying the similar procedure to configuration $\beta$ we obtain a
configuration $\beta'=(B'_1, \dots, B'_n)\in F(\Gamma, n)$ connected to $\beta$
by a continuous collision free motion, such that $B'_i\in {\rm {Int}}(e_0)$ for
all $j=1, \dots, n$.

Next we apply one of the precomputed permutational motions which takes
$\alpha'$ to a configuration $\alpha^{\prime\prime}$ which also lies in the
interior of $e_0$ and has the same order as the configuration $\beta'$.

Finally, the output of the algorithm is the concatenation of (1) the motion
from $\alpha$ to $\alpha'$; (2) the motion from $\alpha'$ to
$\alpha^{\prime\prime}$; (3) the obvious motion from $\alpha^{\prime\prime}$ to
$\beta'$; and (4) the reverse of the motion $\beta\to \beta'$. The motion (3)
exists since both $\alpha^{\prime\prime}$ and $\beta'$ lie on $e_0$ and have
the same order.

The above algorithm is discontinuous: if one of the points $A_j$ is a vertex
then a small perturbation of $A_j$ may lead to a different set of minimal
points and hence to a completely different motion. Note that the vertices of
$T$ which have valence one or two do not cause discontinuity.

Let $S_i\subset F(\Gamma, 2)$ denote the set of all configurations
$\alpha=(A_1, \dots, A_n)$ such that precisely $i$ points among $A_j$ are
essential vertices of $\Gamma$. If we restrict the above algorithm on the pairs
$(A,B)\in S_i\times S_j$ then the result is a continuous function of the input.
The sets $S_i\times S_j$ and $S_{i+1}\times S_{j-1}$ are disjoint and each of
the sets contain no limit points of the other. This follows since the closure
of $S_i$ is contained in the union of the sets $S_j$ with $j\geq i$. We may
define
$$Y_k=\bigcup_{i+j=k} S_i \times S_j, \quad k =0, 1,\dots, 2m(\Gamma)$$
where $m(\Gamma)$ denotes the number of the essential vertices of $\Gamma$. The
algorithm described above is continuous when restricted on each set $Y_k$.
Hence we obtain:

\begin{corollary}
The topological complexity of the algorithm described above is less than or
equal to $2m(\Gamma)+1$ where $m(\Gamma)$ denotes the number of the essential
vertices of $\Gamma$.
\end{corollary}

\section{The Main Result}

Our main result which will be proven later in this paper states:
\begin{theorem}\label{main}
Let $\Gamma$ be a tree having an essential vertex. Let $n$ be an integer satisfying
$n\geq 2m(\Gamma)$ where $m(\Gamma)$ denotes the number of essential vertices of
$\Gamma$. In the case $n=2$ we will additionally assume that the tree $\Gamma$ is not
homeomorphic to the letter $Y$ viewed as a subset of the plane $\mathbf R^2$. Then the
upper bound (\ref{leq}) is exact, i.e.
\begin{eqnarray}\label{equal}
\tc(F(\Gamma, n)) = 2m(\Gamma)+1.
\end{eqnarray}
\end{theorem}

We conjecture that this theorem holds for any connected graph $\Gamma$ having
$\geq 2$ essential vertices without assuming that $\Gamma$ is a tree and for
any number of particles $n\geq 2m(\Gamma)$.

If $\Gamma$ is homeomorphic to the letter $Y$ then $m(\Gamma)=1$ and $F(\Gamma,
2)$ is homotopy equivalent to the circle $S^1$. Hence in this case
$\tc(F(\Gamma,2))=2$, see \cite{F1}. The equality (\ref{equal}) fails in this
case.

From Theorem \ref{thm2} of the next section and Theorem \ref{graph1} it follows
that for any tree $\Gamma$ one has
\begin{eqnarray}\label{tctree}
\tc(F(\Gamma, 2)) = \left\{
\begin{array}{ll}
2, & \, \mbox{if $\Gamma$ is homeomorphic to the letter $Y$},\\ \\
3, &\, \mbox{otherwise.}
\end{array}
\right.
\end{eqnarray}
This example shows that the assumption $n\geq 2m(\Gamma)$ of Theorem \ref{main}
cannot be removed: if $\Gamma$ is a tree with $m(\Gamma)\geq 2$ then
(\ref{equal}) would give $\tc(F(\Gamma,2))=2m(\Gamma)+1\geq 5$ contradicting
(\ref{tctree}).

\begin{figure}[h]
\begin{center}
\resizebox{8cm}{4cm}{\includegraphics[56,336][575,582]{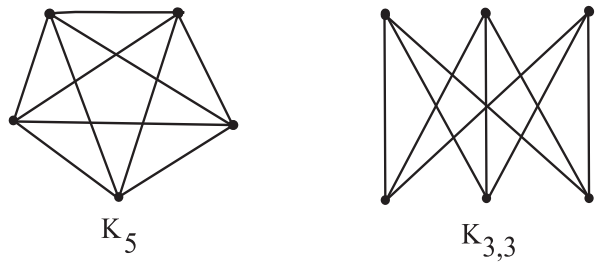}}
\end{center}
\caption{Graphs $K_5$ and $K_{3,3}$.} \label{kkk}
\end{figure}

Here are some more examples. For the graphs $K_5$ and $K_{3,3}$ (see Figure
\ref{kkk}) one has
\begin{eqnarray}
\tc(F(K_5, 2))= \tc(F(K_{3,3},2)) = 5.
\end{eqnarray}
This follows by combining a theorem of A. Abrams \cite{Ab} (stating that the spaces
$F(K_5, 2)$ and $F(K_{3,3},2)$ are homotopy equivalent to closed orientable surfaces of
genus 6 and 4, correspondingly) with Theorem 9 from \cite{F1} computing the invariant
$\tc$ for the surfaces. In these examples the equality(\ref{equal}) is again violated.

\section{Configuration Spaces of Two Particles}

In the present section we explicitly describe the topology of the configuration spaces
$F(\Gamma, 2)$ and $B(\Gamma, 2)=F(\Gamma, 2)/{\mathbf Z}_2$ where $\Gamma$ is a tree.
The latter space is the quotient of $F(\Gamma, 2)$ with respect to the involution
interchanging the particles. The results of this section are used in the proof of Theorem
\ref{main}.

\begin{theorem}\label{thm1} Let $\Gamma$ be a tree having an essential vertex.
Then the space $B(\Gamma, 2)$ is homotopy equivalent to the  wedge of
\begin{eqnarray}\label{sum1}
\frac{1}{2} \sum_{v}(\eta(v)-1)(\eta(v)-2)
\end{eqnarray}
circles where $v$ runs over the vertices of $\Gamma$ and 
the symbol $\eta(v)$ denotes the number of edges incident to the vertex $v$.
\end{theorem}

Note that only essential vertices $v$ contribute nonzero summands to (\ref{sum1}).

In the next theorem we identify the equivariant homotopy type of $F(\Gamma, 2)$
with respect to the canonical involution interchanging the labels of the
particles. Fix a univalent root vertex $u_0$ of $\Gamma$. Then any vertex
$v\not=u_0$ has a well-defined single {\itshape descending} edge which is
incident to it. It is the edge connecting $v$ with the root vertex: removing
the descending edge makes $v$ and $u_0$ lying in different connected
components. The other $\eta(v) -1$ edges incident to the vertex $v$ will be
called {\itshape ascending}.

\begin{theorem}\label{thm2}
Let $\Gamma$ be a tree having an essential vertex. Then the space $F(\Gamma, 2)$ is
homotopy equivalent to the wedge of
\begin{eqnarray}\label{sum2}
\sum_{v}(\eta(v)-1)(\eta(v)-2) -1
\end{eqnarray}
circles where $v$ runs over the vertices of $\Gamma$. Moreover, the ${\mathbf
Z}_2$-equivariant homotopy type of $F(\Gamma, 2)$ can be described as follows. Consider a
1-dimensional cell complex $Y_\Gamma$ having two vertices $A$ and $B$ and a number of
1-dimensional cells connecting $A$ to $B$, each labelled by an ordered pair of distinct
ascending edges of $\Gamma$ incident to an essential vertex $v\in \Gamma$. In total there
are $\sum_v (\eta(v)-1)(\eta(v)-2)$ such 1-cells. The complex $Y_\Gamma$ has an
involution $T:Y_\Gamma\to Y_\Gamma$ which takes $A$ to $B$ and which takes the edge with
the label $(e_i, e_j)$ to the edge with the label $(e_j, e_i)$. Then $F(\Gamma, 2)$ and
$Y_\Gamma$ are ${\mathbf Z}_2$-equivariantly homotopy equivalent.
\end{theorem}

\begin{figure}[h]
\begin{center}
\resizebox{7cm}{2.7cm}{\includegraphics[60, 462][505, 621]{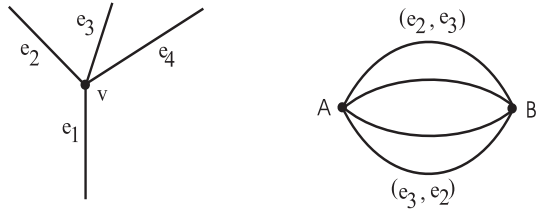}}
\end{center}
\caption{Descending edge $e_1$ and ascending edges $e_2$, $e_3$, $e_4$ incident to a
vertex $v$ (left) and a part of the cell complex $Y_\Gamma$ (right). } \label{descending}
\end{figure}

Figure \ref{descending} illustrates the construction of the cell complex $Y_\Gamma$.

The complex $Y_\Gamma$ can also be described as follows. Consider the standard unit
sphere $S^2\subset \mathbf R^3$ given by the equation $x^2+y^2+z^2=1$. Let $P_\phi$
denotes the 2-dimensional plane spanned by the vectors $(0,0,1)$ and $(\cos\phi,
\sin\phi, 0)$. The intersection $P_\phi\cap S^2$ is a circle $S_\phi^1$ containing the
North and the South poles $(0,0,\pm1)$. Set $\phi_i=\frac{i\pi}{n}$ where $i=0, 1, \dots,
n-1$ and \begin{eqnarray*} n=1/2\cdot \sum_v \eta(v)-1)(\eta(v)-2), \end{eqnarray*} the
sum taken over the vertices of $\Gamma$. Then
\begin{eqnarray}
Y_\Gamma = \bigcup_{i=0}^{n-1} S^1_{\phi_i}.
\end{eqnarray}
The latter space is viewed with the standard antipodal involution.

\section{Sketch of Proof of Theorem \ref{main}.}

First we consider the special case when the tree $\Gamma$ has a single
essential vertex, $m(\Gamma)=1$. Then, by a theorem of R. Ghrist \cite{Gr}, the
configuration space $F(\Gamma, n)$ has homotopy type of a wedge of
\begin{eqnarray}\label{wedge}
1+[(n-1)(\eta-2)-1]\cdot\frac{(n+\eta+2)!}{(\eta-1)!}
\end{eqnarray}
circles where $\eta=\eta(v)$ denotes the number of edges incident to the essential
vertex. Using Theorem \ref{graph1} we find that the topological complexity of the
configuration space $F(\Gamma, n)$ equals either 2 or 3 depending on whether the number
of circles in the wedge is 1 or $> 1$. It is easy to see that (\ref{wedge}) equals 1 if
and only if $\eta=3$ and $n=2$. Since this possibility is excluded by our assumption we
find that $\tc(F(\Gamma, n))=3$ for $m(\Gamma)=1$.

The proof in the case $m(\Gamma)>1$ uses the following lemma:

\begin{lemma}\label{lemma1}
Let $X$ be a topological space and let
\begin{eqnarray}u_1, \dots, u_m, w_1,
\dots, w_m\in H^1(X)\end{eqnarray} be cohomology classes (where
$H^\ast(X)=H^\ast(X;\mathbf k)$ and ${\mathbf k}$ is a field) satisfying
\begin{eqnarray}
u_iw_j =0, \quad i,j=1, \dots, m,
\end{eqnarray}
and such that their cup-products $u_1\dots u_m$ and $w_1\dots w_m\in H^m(X)$
are linearly independent. Then
\begin{eqnarray}
\tc(X)\geq 2m+1.
\end{eqnarray}
\end{lemma}
\begin{proof}
Consider the classes $\bar u_i=1\otimes u_i-u_i\otimes 1$ and $\bar
w_j=1\otimes w_j-w_j\otimes 1$ lying in the tensor product $H^\ast(X)\otimes
H^\ast(X)$. The classes $\bar u_i$ and $\bar w_j$ are zero-divisors. The
product of all these classes in $H^\ast(X)\otimes H^\ast(X)$ equals
\begin{eqnarray}\label{bar}
\prod_{i=1}^m \bar
u_i\cdot \prod_{j=1}^m \bar w_j \, =\,  \pm \prod_{i=1}^m u_i\otimes
\prod_{j=1}^m w_j \, \pm\,  \prod_{j=1}^m w_j \otimes \prod_{i=1}^m u_i.
\end{eqnarray}
Since the classes $\prod_{i=1}^m u_i$ and $\prod_{j=1}^m w_j$ are linearly
independent we may find two linear functionals $\phi, \psi: H^m(X)\to \mathbf
k$ such that $\phi(\prod_{i=1}^m u_i)\not=0$ and $\phi(\prod_{j=1}^m w_j)=0$
while $\psi(\prod_{j=1}^m w_j)\not=0$ and $\psi(\prod_{i=1}^m u_i)=0$. Using
(\ref{bar}) we see that applying $\phi\otimes\psi$ to the product (\ref{bar})
is nonzero which proves the non-triviality of the product (\ref{bar}). Hence we
have a nontrivial product of $2m$ zero divisors. Applying Theorem \ref{lower}
gives the desired result.
\end{proof}

Assume now that $\Gamma$ is a tree having $m(\Gamma)\geq 2$ essential vertices.

Let $v_1, \dots, v_m$ be the essential vertices of $\Gamma$. The space
$F(\Gamma, 2)$ is homotopy equivalent to a wedge of circles (cf. Theorem
\ref{thm2}). For any index $i=1, 2, \dots, m$ we fix a nonzero cohomology class
$\alpha_i\in H^1(F(\Gamma, 2))$ with the only condition that it vanishes on all
the circles in the wedge except the circles associated with the vertex $v_i$,
see formula (\ref{sum2}).

We define two continuous maps
\begin{eqnarray*}
\Phi_i: F(\Gamma, n)\to F(\Gamma, 2), \quad \Psi_i: F(\Gamma, n)\to F(\Gamma,
2),
\end{eqnarray*}
where $i=1, 2, \dots, m$ and
\begin{eqnarray}
\begin{array}{lll}
\Phi_i(x_1, \dots, x_n) &= (x_{2i-1}, x_{2i}),& \\ \\
 \Psi_i(x_1, \dots, x_n) &= \left\{
 \begin{array}{ll}
(x_{2i+1}, x_{2i+2}), &\, \mbox{for}\, i=1, \dots, m-1, \\ \\
(x_1, x_2), &\, \mbox{for}\,
i=m.
\end{array}
\right.
\end{array}
\end{eqnarray}
Finally we denote
\begin{eqnarray}
v_i=\Phi^\ast_i(\alpha_i)\in H^1(F(\Gamma, n)), \quad
w_i=\Psi_i^\ast(\alpha_i)\in H^1(F(\Gamma, n)).
\end{eqnarray}

One checks that the conditions of Lemma \ref{lemma1} are satisfied for the
constructed cohomology classes, i.e. the cup-products $u_1\cdots u_m$ and
$w_1\cdots w_m$ are linearly independent and $u_iw_j=0$. Applying Lemma
\ref{lemma1} completes the proof.

The full details will appear elsewhere.

\section{Conclusion}

The topological invariant $\tc(X)$ imposes important restrictions on the
structure of motion planning algorithms for the mechanical systems having $X$
as their configuration space. $\tc(X)$ bounds from below the {\it order of
instability} of deterministic motion planning algorithms. We prove in this
paper that the number $\tc(X)$ equals the minimal integer $k$ such that there
exists a $k$-valued {\it random motion planning algorithm} in $X$.

The motion planning algorithm in the configuration space $F(\Gamma, n)$ of $n$
distinct points on a tree $\Gamma$ which is described in \S \ref{algorithm} has
the minimal possible topological complexity, as Theorem \ref{main} states. This
algorithm may be used in practical control problems when several objects have
to be moved along a tree $\Gamma$ avoiding collisions.

We observe that for a large number of particles $n$, the topological complexity
of this algorithm depends only on the tree $\Gamma$ and does not depend on the
number of the moving objects $n$.

This result could be compared with the earlier results of \cite{FY} which gives
the topological complexity of the motion planning problem of many objects in
the space $\mathbf R^3$ and on the plane $\mathbf R^2$. The topological
complexity of the motion planning algorithms in these situations depends
linearly on the number of particles; it equals $2n-2$ for $F(\mathbf R^2, n)$
and is $2n-1$ for $F(\mathbf R^3, n)$.

We obtain: {\itshape for a large number of objects which must be simultaneously
controlled avoiding collisions, a great simplification can be achieved by
restricting the motion of the objects to a one-dimensional net.}

This result may potentially have practical applications in some traffic control
problems.

\clearpage
\addcontentsline{toc}{section}{Index}
\flushbottom
\printindex


\begin{thebibliography}{7}
%
\addcontentsline{toc}{section}{References}

\bibitem {Ab} Abrams A. (2002) Configuration spaces of colored graphs.
Geometriae Dedicata {\bf 92}, 185 -- 194


\bibitem{D}  Dold A. (1972) Lectures on Algebraic Topology. Springer - Verlag.


\bibitem{DNF}  Dubrovin B., Novikov S. P. and Fomenko A. T. (1984) Modern Geometry; Methods of
the Homology Theory. Springer-Verlag.



\bibitem{F1} Farber M. (2003) Topological Complexity of Motion Planning.
Discrete and Computational Geometry {\bf29}, 211--221

\bibitem{F2} Farber M. (2004) Instabilities of Robot Motion.
 Topology and its Applications (to appear). Available as preprint cs.RO/0205015.

\bibitem {FTY} Farber M, Tabachnikov S., Yuzvinsky S. (2003)
Topological Robotics: Motion Planning in Projective Spaces. "International
Mathematical Research Notices" {\bf34}, 1853--1870


\bibitem {FY} Farber M., Yuzvinsky S. (2004)
Topological Robotics: Subspace Arrangements and Collision Free Motion Planning.
To appear in AMS volume dedicated to S.P.Novikov's 65th birthday. Available as
preprint math.AT/0210115

\bibitem {Gr} Ghrist R. (2001) Configuration spaces and braid groups on graphs
in robotics. Knots, braids, and mapping class groups -- papers dedicated to
Joan S. Birman, AMS/IP Stud. Adv. Math. {\bf24}, Amer. Math. Soc., Providence,
29 -- 40

\bibitem {GK} Ghrist R., Koditschek D. (2002) Safe cooperative robot dynamics
on graphs. SIAM J. Control Optim. {\bf 40}, 1556 -- 1575

\bibitem {HS} Halperin D., Sharir M. (1995) Arrangements and their applications
in robotics: recent developments. The Algorithmic Foundations of Robotics. K.
Goldberg, D. Halperin, J.C. Latombe and R. Wilson eds., Boston, MA, 495 - 511

\bibitem{L} Latombe J.-C. (1991) Robot Motion Planning. Kluwer Academic Publishers

\bibitem{SS}
 Schwartz J. T., Sharir M. (1983) On the piano movers' problem: II.  General
techniques for computing topological properties of real algebraic manifolds.
Adv. Appl. Math.,  \textbf{4}, 298--351

\bibitem {Sz} Schwarz A.S. (1966) The genus of a fiber space. Amer. Math. Sci. Transl.
{\bf 55}, 49- 140

\bibitem{S} Sharir M. (1997) Algorithmic motion planning. Handbook of Discrete
and Computational Geometry. J.E. Goodman, J. O'Rourke eds. CRC Press, Boca Raton, FL, 733
- 754

\bibitem{Sm} Smale S. (1987) On the topology of algorithms, I.
J. of Complexity,{ \bf 3}, 81-89.

\bibitem {Sp} Spanier E.(1966) Algebraic Topology.



\end{thebibliography}
\end{document}